\documentclass[12pt]{article}
\usepackage{amscd}
\newtheorem{theorem}{Theorem}
\newtheorem{lemma}{Lemma}
\newtheorem{prop}{Proposition}
\newtheorem{cor}{Corollary}

\newtheorem{definition}{Definition}
\newcommand{\rhobar}{\overline{\rho}}

\newcommand{\Qnew}{{Q - {\rm new}}}
\newcommand{\Ad}{{\rm Ad}}
\newcommand{\Ann}{{\rm Ann}}

\usepackage{amssymb}
\usepackage{xy}
\xyoption{all} 
\newif\ifnormalesBeweisEnde

  {\vskip 0.3ex plus 0.5ex minus 0ex \pagebreak[1]
   \global\normalesBeweisEndetrue
   \trivlist
   \item[\hskip\labelsep {\textsc{Proof} \rm of #1}:]}%
  {\ifnormalesBeweisEnde \EndOfBeweis \fi
   \endtrivlist
   \vskip 1ex plus 1ex minus 0ex \pagebreak[2]}
  {\vskip 0.3ex plus 0.5ex minus 0ex \pagebreak[1]
   \global\normalesBeweisEndetrue
   \trivlist
   \item[\hskip\labelsep \textsc{Proof}:]}%
  {\ifnormalesBeweisEnde \EndOfBeweis \fi
   \endtrivlist
   \vskip 1ex plus 1ex minus 0ex \pagebreak[2]}
\def\EndOfBeweis{\hskip .5em \vrule width 1.0ex height 1.0ex depth 0.3ex}

\newcommand{\onto}{\xymatrix @C=1.5pc{\ar@{->>}[r]&}}

\usepackage{latexsym}
\usepackage{graphicx}
\usepackage{psfrag}
\usepackage{amssymb}

\begin{document}

\title{On isomorphisms between deformation rings and Hecke 
rings \\ {\small{\it{Dedicated to the memory of my mother Nalini
B. Khare \\ 30th August 1936--12th March 2002} }}}

\author{Chandrashekhar Khare}

\date{}

\maketitle

\section{Introduction}

It was proven in [W] and [TW] that in many cases
there is an isomorphism between 
certain 
deformation rings and Hecke rings associated to an irreducible, modular, mod 
$p$ representation $$\rhobar:G_{\bf Q} \rightarrow GL_2(k)$$
with $k$ a finite field of characteristic $p \geq 3$ and $G_{\bf Q}$
the absolute Galois group of ${\bf Q}$.
As is well-known, the proven cases were enough to establish the modularity
of semistable elliptic curves.
There were 2 steps in the proof:
\begin{enumerate}
\item Establish an isomorphism between the minimal deformation ring
and minimal Hecke ring (defined below as $R_{\phi}$ and ${\bf T}_{\phi}$):
this was done in [TW] using what are called Taylor-Wiles (TW) systems now.
\item Deduce an isomorphism between non-minimal deformation rings and
Hecke rings using a numerical criterion established in [W] for maps between certain rings to be 
isomorphisms of complete intersections and results about raising levels
of newforms: this was done in [W].
\end{enumerate}

In this article we give a different approach
to establishing  the isomorphisms between deformation rings and Hecke rings
proved in [W] and [TW] in the case when $\rhobar$ is semistable 
that does not use TW systems. While the approach here follows the
over-all strategy of Wiles, instead of the first step above we 
use the isomorphism of certain new quotients ($R_Q^{\Qnew} \simeq {\bf
T}_Q^{\Qnew}$: see Theorem \ref{DT} below) of the deformation rings and
Hecke rings (denoted by $R_Q$ and ${\bf T}_Q$ below) that one would like to prove
isomorphic: this isomorphism results
from the recent Galois cohomological work of R.~Ramakrishna in 
[R] and [R1] and the level raising
results of F.~Diamond and R.~Taylor in [DT] and [DT1]. After this the lifting of the isomorphism
$R_Q^{\Qnew} \simeq {\bf T}_Q^{\Qnew}$ to an isomorphism $R_Q \simeq
{\bf T}_Q$ is carried out by applying the level lowering
results of [RT] (see Theorem \ref{RT} below), and the numerical criterion
alluded to above of Wiles as refined by Lenstra (see Theorem
\ref{numerical criterion} below). In the last section we sketch
another approach that uses a ``level substitution mod $p^n$'' trick
that has some novelty.

\section{Deformation rings and Hecke rings}

Let $\rhobar:G_{\bf Q} \rightarrow 
GL_2(k)$ be continuous and irreducible mod $p$ representation with $k$ a finite
field of characteristic $p>5$ which satisfies the following conditions:

\begin{itemize}
\item    $\rhobar$ is modular, i.e., arises from a newform of
         some level $N$ and weight $k \geq 2$ with respect to an
         embedding $\iota_p:\overline{{\bf Q}} \rightarrow
         \overline{{\bf Q}_p}$ that we fix.

\item ${\rm det}(\rhobar)=\omega$ the mod $p$ cyclotomic     
         character.


\item $$\rhobar|_{G_{\ell}} \simeq \left(\matrix{\omega\chi^{-
1}&*\cr
                                                      0&\chi}\right)$$
for $\chi$ an unramified character of $G_{\ell}$, for $\ell (\neq p) \in S$ with $S$ the
set of places at which $\rhobar$ is ramified (we then have $\chi^2=1$)
and 
with $G_{\ell}$ a decomposition group at $\ell$.

\item $\rhobar|_{G_p}$ is either finite, or not finite but Selmer, 
i.e., of the form $$\left(\matrix{\omega\psi^{-1}&*\cr
                                    0&\psi}\right)$$
for an unramified character $\psi$, with the further condition 
that its restriction to the inertia subgroup 
is not the generic fibre of a finite flat group scheme 
over the ring of integers of the maximal unramified extension 
of ${\bf Q}_p$ (we then have $\psi^2=1$).
\end{itemize}

Note that this implies that $\Ad^0(\rhobar)$ is absolutely irreducible
and the order of ${\rm im}(\rhobar)$ is divisible by $p$ (otherwise $\rhobar$ has Serre weight 2
and Artin conductor 1 by assumptions on determinant and ramification above) 
which give the properties of $\rhobar$ required in [R] and [R1] (we
thank the referee for pointing this out).

\subsection{Deformation rings}

We briefly recall the existence of 
certain deformation rings parametrising lifts of $\rhobar$ with
given local conditions, referring to [Ma] and [KR] for more details. We denote
by $\varepsilon$ the $p$-adic cyclotomic character.
For a complete Noetherian local $W(k)$-algebra $R$ with residue field $k$, a representation
$\rho:G_{\bf Q} \rightarrow GL_2(R)$ that reduces modulo the 
maximal ideal to give $\rhobar$, is said to be minimally ramified at
at a prime $\ell$ if:
\begin{enumerate} \item $\ell$ is not in $S$,
then $\rho$ is unramified at $\ell$ \item $\ell \in S$ and
$\ell \neq p$, then if we denote by $I_{\ell}$ an inertia subgroup at
$\ell$, 
$\rho|_{I_{\ell}}$ is isomorphic to  
$$\left(\matrix{1&*\cr
                       0&1}\right),$$ 
\item if $\ell=p$ then $\rho|_{I_p}$  is finite if 
$\rhobar$ is finite, and is otherwise of the form
$$\left(\matrix{\varepsilon&*\cr
                0&1}\right).$$ 
\end{enumerate}

For a finite set of primes $Q$, not congruent to $\pm 1$ mod 
$p$, at which $\rhobar$ is unramified, and $\rhobar|_{G_q}$ 
is of the form $${\pm}\left(\matrix{\varepsilon&0\cr
                                    0&1}\right),$$
consider the ring $R_Q^{\Qnew}$ (that is called $R_{S \cup Q}^{\Qnew}$
in [KR]) that is
the universal ring that parametrises (isomorphism classes of) 
lifts of $\rhobar$ 
to $GL_2(R)$ for complete Noetherian local $W(k)$-algebras $R$ 
with residue field $k$, such that these are

\begin{itemize}

\item Unramified outside $S \cup Q$

\item Minimally ramified, in the above sense, at primes in the set $S$ 
that consists of the ramified primes of $\rhobar$ 

\item At primes $q \in Q$ on restriction to $G_q$ have the form

$${\pm} \left(\matrix{\varepsilon&*\cr
                      0&1}\right)$$

\item With determinant $\varepsilon$, the $p$-adic cyclotomic character. 

\end{itemize}

Such a universal ring $R_Q^{\Qnew}$, that is a complete Noetherian
local
$W(k)$-algebra with residue field $k$, exists as the conditions above 
are conditions in the sense of Section 23 of [Ma] because
of our assumption that the primes in $Q$ are not congruent to
1 mod $p$ (see Lemma 2 of [KR] for more details), and we denote
the corresponding universal representation $G_{\bf Q} \rightarrow
GL_2(R_Q^{\Qnew})$
by $\rho_Q^{\Qnew}$. If we delete
the third condition on primes in $Q$ then we denote the corresponding
ring $R_Q$. We have the natural surjective map $R_Q \rightarrow
R_Q^{\Qnew}$ of local $W(k)$-algebras. More generally for a subset $\alpha$ of $Q$ if we impose
the third condition only for primes in $\alpha$
then we denote the corresponding deformation ring by
$R_Q^{\alpha-{\rm new}}$
and the corresponding universal representation $G_{\bf Q} \rightarrow
GL_2(R_Q^{\alpha-{\rm new}})$ by $\rho_Q^{\alpha-{\rm new}}$:
we have surjective maps $R_Q \rightarrow
R_Q^{\alpha-{\rm new}} \rightarrow R_Q^{\Qnew}$ of local $W(k)$-algebras.
If $Q$ is empty we denote the corresponding ring by $R_{\phi}$ and
call it the {\it minimal deformation ring}, and denote by
$\rho_{\phi}$ the corresponding universal representation.

We denote by $G_{S \cup Q}$ the Galois
group of the maximal extension of $\bf Q$ in $\overline{\bf Q}$
unramified outside $S$ and $Q$ and the infinite place.

\begin{definition}

The Selmer group $H^1_Q(G_{S \cup Q},\Ad^0(\rhobar))$ is 
the kernel of the restriction map $$H^1(G_{S \cup 
Q},\Ad^0(\rho)) \rightarrow \oplus_{v \in S \cup 
Q}H^1(G_v,\Ad^0(\rhobar))/{\cal N}_v.$$ 

\end{definition}

Here for $v \neq p$ and 
$v \in S$, ${\cal N}_v$ is described as the image 
under the inflation map of 
$H^1(G_v/I_v,\Ad^0(\rhobar)/{\Ad^0(\rhobar)}^{I_v})$.
For $v=p$ we define it to be either $H^1_{Se}(G_p,\Ad^0(\rhobar))$
or $H^1_{fl}(G_p,\Ad^0(\rhobar))$, according to whether
$\rhobar|_{I_p}$ is not finite or it is finite, using the notation of Section 4.1
of [deS]. 
For $v \in Q$, ${\cal N}_v$ is described as the 
subspace of $H^1(G_v,\Ad^0(\rhobar))$ generated by the 
cocycle that (in a suitable choice of basis of $V_{\rhobar}$, the
2-dimensional $k$-vector space that affords $\rhobar$, and viewing
$\Ad^0(\rhobar)$ as a subspace of ${\rm End}(V_{\rhobar})$) 
sends $\sigma \rightarrow 0$ and $\tau$ to
$$\left(\matrix{0&1\cr
                0&0}\right)$$
where $\sigma$ and $\tau$ generate the tame quotient of 
$G_q$, and satisfy the relation $\sigma\tau\sigma^{-
1}=\tau^q$ (this is the ``null cocycle'' 
denoted by $s_q$ in Section 3 of [R]).

\begin{prop}\label{ts}
  If $\sf m_Q$ is the maximal ideal of $R_Q^{Q-new}$, then we have 
  an isomorphism $$H^1_Q(G_{S \cup Q},\Ad^0(\rhobar))   
  \simeq {\rm Hom}({\sf m}_Q/({\sf m}_Q^2,p),k).$$
\end{prop}  

\noindent{\bf Proof.} The proof is routine and follows from 
the standard identification of the cotangent space of 
deformation rings with the corresponding Selmer groups
(see Theorem 15 of [deS] for instance or [Ma1] in [FLT]).

\subsection{Hecke rings}

Let $N(\rhobar)$ be the prime to 
$p$ part of the Artin conductor of $\rhobar$, and let $\delta=0$ 
or 1 according to whether $\rhobar$ is finite or not at $p$.
Consider a finite set
of primes $Q=\{q_2,\cdots ,q_{2r}\}$ of odd cardinality (the reason
for this eccentric numbering will be apparent) such that $\rhobar$ is
unramified
at primes in $Q$, $q_i$ is not $\pm 1$ mod $p$ for $q_i \in Q$  
and such that ${\rm tr}(\rhobar({\rm Frob}_{q_i}))=
\pm (q_i+1)$ for $q_i \in Q$ with ${\rm Frob}_{q_i}$ a (arithmetic) Frobenius element at $q_i$.
For the subset of primes of even cardinality
$\{{q_1,\cdots,q_{2(r-s)}}\}$ (for any $0 \leq s \leq r$)
of $Q \cup \{q_1\}$, with $q_1=t$ a fixed prime 
dividing $N(\rhobar)p^{\delta}$ ($>1$) that we choose once and for all,
let $B_{q_1,\cdots,q_{2(r-s)}}$
be the indefinite quaternion algebra
ramified exactly at the set $\{q_1,\cdots,q_{2(r-s)}\}$. 
For the congruence subgroup $\Gamma_0(N(\rhobar)p^{\delta}q_{2(r-s)+1} \cdots
q_{2r}/t)$ of $B_{q_1,\cdots,q_{2(r-s)}}$ denote
the Jacobian of the corresponding Shimura curve
by $J^{q_1,\cdots,q_{2(r-s)}}$. This has an action of the Hecke
operators $T_r$ for primes $r$ coprime 
to the primes in $Q$ and coprime to $N(\rhobar)p$.
By virtue of the main theorem of [DT] and the Jacquet-Langlands
correspondence there is a maximal ideal $\sf m$ of the $W(k)$-algebra generated
by the induced action of these $T_r$'s on 
${\rm Ta}_p(J^{q_1,\cdots,q_{2(r-s)}})
\otimes_{{\bf Z}_p}W(k)$ (where ${\rm Ta}_p(J^{q_1,\cdots,q_{2(r-s)}})$
is as usual the $p$-adic Tate module of $J^{q_1,\cdots,q_{2(r-s)}}$) that gives rise to the
representation $\rhobar$, i.e., $T_r -\widetilde{{\rm tr}(\rhobar({\rm
Frob}_r))} \in \sf m$ for $\widetilde{{\rm tr}(\rhobar({\rm
Frob}_r))}$ an arbitrary lift to $W(k)$ of ${\rm tr}(\rhobar({\rm
Frob}_r))$ and ${\rm Frob}_r$ a (arithmetic) Frobenius element at $r$. 
We denote by ${\bf T}_Q^{\{q_1,\cdots,q_{2(r-s)}\}-{\rm new}}$ the 
completion of this Hecke algebra at $\sf m$: it is finite flat as
a $W(k)$ module. 

To define the {\it minimal Hecke ring} ${\bf T}_{\phi}$, we consider  
${\rm Ta}_p(J_0(N(\rhobar)p^{\delta})) \otimes_{{\bf Z}_p} W(k)$ where
${\rm Ta}_p(J_0(N(\rhobar)p^{\delta}))$ is the $p$-adic Tate module
of the Jacobian $J_0(N(\rhobar)p^{\delta})$ of the modular curve
$X_0(N(\rhobar)p^{\delta})$: this has an action of 
the Hecke operators $T_r$ for primes $r$ that are coprime
to $N(\rhobar)p$ and the primes in the chosen finite set of primes
$Q$. The $W(k)$-algebra generated by the 
images of these $T_r$'s in the $W(k)$-endomorphisms of ${\rm Ta}_p(J_0(N(\rhobar)p^{\delta})) \otimes_{{\bf Z}_p} W(k)$
is finite and flat over $W(k)$. By results of
[Ri] (see [Ri1] for a proof that does not invoke multiplicity one
results) there is a maximal ideal $\sf m$ of this Hecke algebra that gives
rise to $\rhobar$ in the sense explained above. We denote by 
${\bf T}_{\phi}$ the completion at $\sf m$ of this Hecke algebra.

\begin{lemma}
  We have representations
  $\rho_{Q,{\rm mod}}^{\{q_1,\cdots,q_{2(r-s)}\}-new}:G_{\bf Q}
  \rightarrow GL_2({\bf T}_Q^{\{q_1,\cdots,q_{2(r-s)}\}-{\rm new}})$
  (resp., $\rho_{\phi,{\rm mod}}:G_{\bf Q} \rightarrow GL_2({\bf T}_{\phi})$)
  that are unramified outside primes in $S \cup Q$ and the infinite
  place, and are characterised by ${\rm
  tr}(\rho_{Q,{\rm mod}}^{\{q_1,\cdots,q_{2(r-s)}\}-new}({\rm Frob}_r))=T_r$
  (resp., ${\rm tr}(\rho_{\phi,{\rm mod}}({\rm Frob}_r))=T_r$), that arise by
  uniquely determined specialisations of the universal representations
  $\rho_Q^{ \{ q_2,\cdots, q_{2(r-s)}\}-{\rm new} }:G_{\bf Q} \rightarrow 
  GL_2(R_Q^{\{q_2,\cdots,q_{2(r-s)}\}-{\rm new}})$ (resp.,
  $\rho_{\phi}:G_{\bf Q} \rightarrow GL_2(R_{\phi})$).
  The aforementioned specialisation maps
  $R_Q^{\{q_2,\cdots,q_{2(r-s)}\}-{\rm new}}
  \rightarrow {\bf T}_Q^{\{q_1,\cdots,q_{2(r-s)}\}-{\rm new}}$ (resp.,
  $R_{\phi} \rightarrow {\bf T}_{\phi}$) of local $W(k)$-algebras are surjective.
\end{lemma}

\noindent{\bf Proof:} The existence of the representations
  $\rho_{Q,{\rm mod}}^{\{q_1,\cdots,q_{2(r-s)}\}-{\rm new}}$ and
  $\rho_{\phi,{\rm mod}}$ with the asserted properties
  follows from the results of [Ca] and the rest from
  the universal properties of deformation rings and the properties of
  the representations $\rho_{Q,{\rm mod}}^{\{q_1,\cdots,q_{2(r-s)}\}-{\rm new}}$ and
  $\rho_{\phi,{\rm mod}}$ deduced from [Ca1]: for instance 
  $\rho_{Q,{\rm mod}}^{\{q_1,\cdots,q_{2(r-s)}\}-{\rm new}}|_{G_q}$ is
  of the form $${\pm} \left(\matrix{\varepsilon&*\cr
                      0&1}\right)$$ by [Ca1]. 

\subsection{Isomorphisms between deformation rings and Hecke rings}

\begin{definition}\label{special}
  For $\rho:G_{\bf Q} \rightarrow GL_2(W(k))$ or 
  $\rho:G_{\bf Q} \rightarrow GL_2(W(k)/p^n)$ we say that
  $\rho$ is {\bf special} at a prime $q$ if $q$
  is not congruent to $\pm 1$ mod $p$, and $$\rho|_{G_q} \simeq
  {\pm} \left(\matrix{\varepsilon&*\cr
                0&1}\right).$$
\end{definition}

The following theorem and its proof we owe to Ravi Ramakrishna:
we use the notation of Sections 2.1 and 2.2.

\begin{theorem}\label{DT}
  There exists a finite set of primes $Q=\{q_2,\cdots,q_{2r}\}$ 
  of odd cardinality for some integer $r$ such that for $q \in
  Q$, $\rhobar$ is unramified and special at $q$ and 
  we have the sequence of isomorphisms $R_Q^{\Qnew} \simeq {\bf
  T}_Q^{\Qnew} \simeq W(k)$ of local $W(k)$-algebras.
\end{theorem}

\noindent{\bf Proof.} Using Proposition 21 of [KR] we get 
a finite set of primes $Q=\{q_2,\cdots,q_{2r}\}$ such that for $q \in Q$, $q$
is special for $\rhobar$ and such that 
the map
$$H^1(G_{S \cup Q},\Ad^0(\rhobar)) \rightarrow \oplus_{v \in S \cup Q}
H^1(G_v,\Ad^0(\rhobar))/{\cal N}_v$$
is an isomorphism. Thus by Proposition \ref{ts},
$R_Q^{\Qnew}$ is a quotient of $W(k)$.
The main theorem of [DT] 
yields a $p$-adic modular lift of $\rhobar$, that
arises from a specialisation of the universal representation 
$\rho_Q^{\Qnew}:G_{\bf Q} \rightarrow GL_2(R_Q^{\Qnew})$. As by Lemma
1 we {\it a priori} do have a 
surjection $R_Q^{\Qnew} \rightarrow {\bf T}_Q^{\Qnew}$ the
isomorphisms
in the statement follow.

\begin{definition}\label{ramaux}
  A finite set of primes $Q$ such that $R_Q^{\Qnew} \simeq 
  W(k)$ as in the theorem is called an auxiliary set and the
  corresponding representation $G_{\bf Q} \rightarrow GL_2(W(k))$
  is denoted by $\rho_Q^{\Qnew}$.
\end{definition}

A new proof of the following theorem of [W] and [TW], that will
deduce it from Theorem \ref{DT}, 
will be given in the next section  
and is the main contribution of this article.

\begin{theorem}\label{main}
  We have an isomorphism $R_Q \simeq {\bf T}_Q$ of complete
  intersection rings.
\end{theorem}

We have the following corollary proved in [W] and [TW] that is deduced from
Theorem \ref{main} using the methods of Section 2 of
[W], i.e., level raising methods as in [Ri2] and [DT] and the
numerical isomorphism criterion of Wiles.

\begin{cor}\label{modularity}
  Let $\rho:G_{\bf Q} \rightarrow GL_2({\cal O})$ be a lift of 
  $\rhobar$, with $\cal O$ the ring of integers of a finite extension
  of ${\bf Q}_p$, such that $\rho$ is 
  unramified at almost all places, with determinant
  the $p$-adic cyclotomic character $\varepsilon$, such that $\rho|_{I_p}$
  either arises from a $p$-divisible group scheme
  or $\rho|_{I_p}$ is of the form
  $$\left(\matrix{\varepsilon&*\cr
                  0&1}\right).$$                
  Then $\rho$ is modular.
\end{cor}

We also have the following corollary that can be deduced from
Theorem \ref{main}: this deduction is in the appendix by Gebhard
B\"ockle [B] to this paper.

\begin{cor}
   We have an isomorphism $R_{\phi} \simeq {\bf T}_{\phi}$ of complete
   intersection rings.
\end{cor}

\section{Proof of Theorem \ref{main}}

\subsection{Some preliminaries}

We begin by stating 2 results needed below.
The first is the vital commutative algebra result
of [W] as refined by Lenstra (cf. Theorem 5.3 of [DDT]).
(We say that a $W(k)$-algebra $R$
that is a complete Noetherian local ring,
and is finite flat as a $W(k)$ module, with residue field $k$,
is a complete intersection
if $R \simeq W(k)[[T_1,\cdots,T_n]]/(f_1,\cdots,f_n)$ for some $n$.)

\begin{theorem}\label{numerical criterion}
  
  Let $R$ and $T$ be complete Noetherian local $W(k)$-algebras with residue
  field $k$, $T$ finite flat as a $W(k)$ module, 
  and $\phi:R \rightarrow T$ a surjective map of local $W(k)$-algebras. Let
  $\pi:T \rightarrow W(k)$ be a homomorphism of local $W(k)$-algebras, 
  and set $\Phi_R={\rm ker}(\pi\phi)/{\rm ker}(\pi\phi)^2$ and 
  $\eta_T=\pi(\Ann_T({\rm ker}(\pi)))$.
  Then $|W(k)/\eta_T| \leq |\Phi_R|$, where in the case
  $|W(k)/\eta_T|$ is infinite this is interpreted to mean
  that $\Phi_R$ is infinite too. Assume that $\eta_T$ is not
  zero. Then the following are equivalent:

\begin{itemize}

\item The equality $|\Phi_R|=|W(k)/\eta_T|$ is satisfied.

\item The rings $R$ and $T$ are complete intersections, and $\phi$ is
      an isomorphism.

\end{itemize}

\end{theorem}

The second result we need is an upper bound on the orders of the tangent
spaces at $\pi$ of the deformation rings $R_Q^{\alpha-{\rm new}}$ where $\alpha$ is
a subset of the finite set $Q$, $Q$ is as in Theorem \ref{DT} and
$\pi$ is the map that is obtained by composing the natural map
$R_Q^{\alpha-{\rm new}}
\rightarrow R_Q^{\Qnew}$ with the isomorphism 
$R_Q^{\Qnew} \simeq W(k)$ of Theorem \ref{DT}. We define
$\Phi_{\pi}(R_{Q}^{\alpha-{\rm new}})={\rm ker}(\pi)/{\rm ker}(\pi)^2$.
The representation $\rho_Q^{Q-new}:G_{\bf Q} \rightarrow GL_2(W(k))$
is residually absolutely irreducible and thus 
is integrally well-defined by [Ca]. Looking at 
the semistable representation $\rho_Q^{\Qnew}|_{G_q}$ for any prime $q \in Q$ 
we have ideals $(x_{q})$ of $W(k)$, 
where $\rho_Q^{\Qnew}(t_{q})$, with $t_{q}$ generating
the unique ${\bf Z}_p$ quotient of $I_{q}$, is of the form
$$\left(\matrix{1&x_{q}\cr
                0&1\cr}\right),$$ for some $x_q \in W(k)$ and $(x_q)$ does not
depend on the choice of an integral model for $\rho_Q^{Q-new}$. We
note that $x_q \neq 0$ (as $R_Q^{\Qnew} \simeq {\bf T}_Q^{\Qnew}$ by 
Theorem \ref {DT} and then use [Ca1]).
We have the following result that is proved 
in Section 5 of [KR] by computations of orders of some local
Galois cohomology groups.

\begin{prop}\label{tangent space}
  $|\Phi_{\pi}(R_{Q}^{\alpha-new})| \leq \Pi_{q \in Q \backslash \alpha}
  |W(k)/(x_{q})|$.
\end{prop}

The strategy of the proof of Theorem \ref{main} is to deduce it from Theorem \ref{DT} using
the above 2 results and the result of [RT] to compute change of
$\eta$-invariants when we relax newness conditions 
on the deformation and Hecke rings two primes at a time.

We use the notation of Section 2 and consider an auxiliary set
$Q=\{q_2,\cdots,q_{2r}\}$
provided by Theorem \ref{DT}.
Then the isomorphism $\pi:{\bf T}_Q^{\Qnew} \rightarrow W(k)$ of Theorem \ref{DT}
induces a morphism ${\bf T}_Q^{\{q_1,\cdots,q_{2(r-s)}\}-{\rm new}}
\rightarrow W(k)$ of local $W(k)$-algebras using the Jacquet-Langlands
correspondence: by an abuse of notation we
again denote this morphism by $\pi$.

Let $\xi_s:J_{\Gamma}^{q_1,\cdots,q_{2r-2s}} \rightarrow A_s$ 
be the optimal quotient that corresponds to the 
morphism  $\pi$ ($0 \leq s \leq r-1$): thus $A_s$ is an abelian
variety
over ${\bf Q}$ such that $\rho_Q^{\Qnew}$ arises from its $p$-adic
Tate module, the map $\xi_s$ is equivariant for the Hecke action
and its kernel is connected. 
Analogously we have the optimal
quotient $\xi_s:J_{\Gamma}^{q_1,\cdots,q_{2r-2(s+1)}} 
\rightarrow A_{s+1}$. (The optimal quotients that we consider
are all isogenous to each other by Faltings' isogeny theorem, 
via Hecke-equivariant isogenies defined over ${\bf Q}$.) We
need to study the maps ${\rm Ta}_p(J^{q_1,\cdots,q_{2r-2s}}) 
\rightarrow {\rm Ta}_p(A_s)$
and ${\rm Ta}_p(J^{q_1,\cdots,q_{2r-2(s+1)}})
\rightarrow {\rm Ta}_p(A_{s+1})$ that $\xi_s$ and $\xi_{s+1}$ induce on 
the corresponding Tate
modules. With the help of these we will be able to 
compare the $\eta$-invariant of the homomorphism
$\pi: {\bf T}_Q^{\{q_1,\cdots,q_{2r-2(s+1)}\}-{\rm new}} \rightarrow
W(k)$ (corresponding to $\xi_{s+1}$) to the $\eta$-invariant of the homomorphism
$\pi: {\bf T}_Q^{\{q_1,\cdots,q_{2r-2s}\}-{\rm new}} \rightarrow
W(k)$ (corresponding to $\xi_s$). 
To ease the notation, we set $J:=J^{q_1,\cdots,q_{2r-2s}}$
and $J':=J^{q_1,\cdots,q_{2r-2(s+1)}}$, $A_s:=A$ and $A_{s+1}:=A'$,
$\xi:=\xi_s$
and $\xi':=\xi_{s+1}$ and the
corresponding Hecke algebras
${\bf T}_Q^{\{q_1,\cdots,q_{2r-2s}\}-{\rm new}}$ 
and ${\bf T}_Q^{\{q_1,\cdots,q_{2(r-2(s+1))}\}-{\rm new}}$ by $\bf T$ 
and ${\bf T}'$ respectively. We also set $q_{2(r-s)}=q$ and $q_{2(r-s)-1}=q'$.

We have induced maps $\xi'^*:{\rm Ta}_p({(A')}^d)_{{\sf m},W(k)}
\rightarrow {\rm Ta}_p(J')_{{\sf m},W(k)}$,
and $\xi'_*:{\rm Ta}_p(J')_{{\sf m},W(k)} \rightarrow {\rm Ta}_p(A')_{{{\sf m}},W(k)}$ (injective
with torsion-free cokernel and surjective respectively by optimality, and their analogs for
$\xi$), where for instance ${\rm Ta}_p({(A')}^d)_{{\sf m},W(k)}:=({\rm
Ta}_p({(A')}^d) \otimes_{{\bf
Z}_p} W(k))_{\sf m}$ is the localisation at $\sf m$ of $p$-adic Tate module 
of the dual abelian variety $(A')^d:={\rm Pic}^0(A')$ of $A'$ tensored 
with $W(k)$, and the other symbols are defined analogously.
(It's worth noting that if $k$ is the field
of definition of $\bar{\rho}$, $\cal O$ is the
image of the Hecke algebra in ${\rm End}(A')$ and $\lambda$ is the
prime of $\cal O$ corresponding to $\sf m$, then ${\cal O}_\lambda \cong W(k)$,
${\rm Ta}_p(A')_{{\sf m},W(k)} \cong {\rm Ta}_p(A')_\lambda$, etc.) Observe
that ${\rm Ta}_p((A')^d)_{{{\sf
m}},W(k)}$ 
is free of rank 2 over the quotient of ${\bf T}'$ cut out by its
action on it: this quotient is isomorphic to $W(k)$
(a section of the structure map $W(k) \rightarrow {\bf T}'$). The module
${\rm Ta}_p((A')^d)_{{{\sf m}},W(k)}$ is
(non-canonically) isomorphic to ${\rm Ta}_p(A')_{{{\sf m}},W(k)}$ 
as a $W(k)[G_{\bf Q}]$ module by
[Ca]. We identify these by choosing any such isomorphism: the
ambiguity is only up to elements of $W(k)^*$ which is immaterial in the
calculations below. 
Analogously we have the maps induced by $\xi$ 
involving a choice of a $W(k)[G_{\bf Q}]$-isomorphism between localisations at $\sf m$ of 
the Tate modules of $A$ and its dual. 
The maps $\xi'_*\xi'^*$ and  $\xi_*\xi^*$
of ${\rm Ta}_p(A')_{{{\sf m}},W(k)}$ and ${\rm Ta}_p(A)_{{\sf m},W(k)}$ commute with the
$W(k)[G_{\bf Q}]$-action 
and by the irreducibility of this action
can be regarded as given by multiplication by elements of $W(k)$. We denote the
corresponding ideals of $W(k)$ by $(\xi'_*\xi'^*)$ and  $(\xi_*\xi^*)$.

\subsection{A result of Ribet and Takahashi}

The following theorem which follows from an adaptation of
the methods of [RT] to our situation is key to the deduction of 
Theorem \ref{main} from Theorem \ref{DT}.

\begin{theorem}\label{RT}
  $(\xi'_*\xi'^*)=(x_{q}x_{q'})(\xi_*\xi^*)$ (where when $s=r-1$ we
  declare $x_{q'}$ to be a unit).
\end{theorem} 

\noindent{\bf Proof:} Our argument is
essentially that of Sections 2 and 3 of [RT], and below we just
make the necessary adaptations of their arguments, that
are in the case of optimal elliptic curve quotients of the above
Jacobians, to make them go through in the case of optimal abelian
variety quotients that we consider here
which do not come equipped with canonical principal polarisations. 
The key point in [RT] as well as here is 
to consider character groups of the toric
parts of the reductions
of $J$ and $J'$ at $q'$ and $q$ respectively that are related by
Ribet's exact sequence (\ref{ribet}) below. The general facts about
semistable abelian varieties that we use are found in [G].

Consider ${\cal X}(J',q)_{{\sf m},W(k)}={({\cal X}(J',q) \otimes_{{\bf
Z}} W(k))}_{\sf m}$ and 
${\cal X}(J,q')_{{\sf m},W(k)}={({\cal X}(J,q') \otimes_{{\bf Z}}
W(k))}_{\sf m}$,
where ${\cal X}(J',q)$ and ${\cal X}(J,q')$ are the character
groups of the tori associated to the special fibre of  N\'eron models
of the abelian varieties
$J'$ and $J$, that are defined over ${\bf Q}$, at $q$ and $q'$ respectively. 
These have functorial actions
by the Hecke algebras ${\bf T}'$ and ${\bf T}$ respectively.
We analogously consider the character
groups (tensored with $W(k)$ and localised at $\sf m$) ${\cal X}(A,q')_{{\sf m},W(k)}$ and 
${\cal X}(A',q)_{{\sf m},W(k)}$ of
the tori of the reduction of the N\'eron model of the abelian
varieties $A,A'$
at the primes $q'$ and $q$ respectively. Note that as $A'$ has purely
toric reduction at $q$, we have the exact sequence of $W(k)[G_q]$-modules
$$ 0 \rightarrow {\rm Hom}_{W(k)}({\cal X}(A',q)_{{\sf m},W(k)},W(k)(1))
\rightarrow {\rm Ta}_p(A')_{{\sf m},W(k)} \rightarrow {\cal
X}((A')^d,q)_{{\sf m},W(k)}  \rightarrow 0$$
and the analogous exact sequence for the Tate module of the dual abelian variety
$(A')^d$. Now we claim that
the $W(k)[G_{\bf Q}]$-equivariant
isomorphism between ${\rm Ta}_p(A')_{{\sf m},W(k)}$ and
${\rm Ta}_p((A')^d)_{{\sf m},W(k)}$ chosen above induces yields a commutative diagram 
\begin{equation}\label{duality}
\begin{CD}
0 @> >> {\rm Hom}_{W(k)}({\cal X}(A',q)_{{\sf m},W(k)},W(k)(1)) @> >>
 {\rm Ta}_p(A')_{{\sf m},W(k)}
 @> >> {\cal X}((A')^d,q)_{{\sf m},W(k)}  @>  >> 0
\\ 
&& @V VV @V VV @V VV 
\\
0 @> >>{\rm Hom}_{W(k)}({\cal X}((A')^d,q)_{{\sf m},W(k)},W(k)(1)) @>
 >> {\rm Ta}_p((A')^d)_{{\sf m},W(k)}
 @> >>{\cal X}(A',q)_{{\sf m},W(k)} @>  >> 0
\end{CD}
\end{equation} 
in which all the vertical arrows are ismorphisms.
This follows from the fact that the middle vertical arrow
being $W(k)[G_q]$-equivariant, it induces an isomorphism $${\rm
Hom}_{W(k)}({\cal X}(A',q)_{{\sf m},W(k)},W(k)(1)) \otimes {\bf Q}_p \simeq 
{\rm Hom}_{W(k)}({\cal X}((A')^d,q)_{{\sf m},W(k)},W(k)(1)) \otimes
{\bf Q}_p.$$ Noting 
that the cokernels of the maps
in the horizontal exact sequences above 
are torsion-free establishes the claim. (When $q \neq p$ this follows
more simply from noting that the first term in each sequence is the
module of $I_q$-invariants.)

We have an injective map 
$\xi_X'^*:{\cal X}(A',q)_{{\sf m},W(k)} \rightarrow 
{\cal X}(J',q)_{{\sf m},W(k)}$ induced by pull-back from $\xi'$ and similarly a 
surjective map $\xi'_{X,*}:{\cal X}(J',q)_{{\sf m},W(k)}
\rightarrow {\cal X}((A')^d,q)_{{\sf m},W(k)}$ induced by pull-back from the map dual
to $\xi'$ from $(A')^d \rightarrow J'$. Note that as $J'$ 
has semistable reduction at $q$, by results of [G] we have an
inclusion
${\rm Hom}_{W(k)}({\cal X}(J',q)_{{\sf m}, W(k)},W(k)(1))
\rightarrow {\rm Ta}_p(J')_{{\sf m},W(k)}$ that is ${\bf T}'[G_q]$-equivariant. Further, the map
$\xi'^*:{\rm Ta}_p((A')^d)_{{\sf m},W(k)} \rightarrow {\rm Ta}_p(J')_{{\sf m},W(k)}$ (that is induced from
the map $(A')^d \rightarrow J'$ dual to $\xi'$ as above) induces a
Hecke equivariant map
$${\rm Hom}_{W(k)}({\cal X}((A')^d,q)_{{\sf m},W(k)},W(k)(1)) \rightarrow {\rm
Hom}_{W(k)}({\cal X}(J',q)_{{\sf m}, W(k)},W(k)(1))$$ 
of $W(k)[G_q]$-modules
that is ${\rm Hom}(\xi'_{X,*},W(k)(1))$. 
These considerations show that the sequence of maps
$${\rm Ta}_p(A')_{{\sf m},W(k)} \simeq {\rm Ta}_p((A')^d)_{{\sf
m},W(k)} \rightarrow  {\rm Ta}_p(J')_{{\sf m},W(k)} \rightarrow 
{\rm Ta}_p(A')_{{\sf m},W(k)},$$ (where the first map is the
isomorphism we have fixed) induces the sequence of maps
$${\rm Hom}_{W(k)}({\cal X}(A',q)_{{\sf m},W(k)},W(k)(1))
\simeq {\rm Hom}_{W(k)}({\cal X}((A')^d,q)_{{\sf m},W(k)},W(k)(1))$$ $$\rightarrow {\rm
Hom}_{W(k)}({\cal X}(J',q)_{{\sf m}, W(k)},W(k)(1)) \rightarrow
{\rm Hom}_{W(k)}({\cal X}(A',q)_{{\sf m},W(k)},W(k)(1))$$ 
(where for instance the first map is the
isomorphism deduced from (\ref{duality}) and the last map is 
${\rm Hom}(\xi_X'^*,W(k)(1))$). As the composite of the
maps in the first sequence is multiplication by $\xi'_*\xi'^* \in
W(k)$, we deduce that the composite of the
maps in the second sequence is also multiplication by $\xi'_*\xi'^*$.

The isomorphism of $W(k)[G_q]$-modules ${\rm Hom}_{W(k)}({\cal X}((A')^d,q)_{{\sf
m},W(k)},W(k)(1)) \simeq 
{\rm Hom}_{W(k)}({\cal X}(A',q)_{{\sf m},W(k)},W(k)(1))$
induces an isomorphism 
${\cal X}(A',q)_{{\sf m},W(k)} \simeq {\cal X}((A')^d,q)_{{\sf m},W(k)}$
which we use to identify these two latter modules. 
Using this isomorphism we will regard 
$\xi'_{X,*}$ as a map $\xi'_{X,*}:{\cal X}(J',q)_{{\sf m},W(k)}
\rightarrow {\cal X}(A',q)_{{\sf m},W(k)}$. 
From what we have said above it follows that
the composition $\xi'_{X,*}\xi_X'^*$ acts on the
rank one $W(k)$-module ${\cal X}(A',q)_{{\sf m},W(k)}$ by 
multiplication by an element of $W(k)$ which generates
the ideal that we denoted by $(\xi'_*\xi'^*)$ above. From this point onwards we
will abuse notation and denote the maps $\xi_X'^*$ and $\xi'_{X,*}$
by $\xi'^*$ and $\xi'_*$ respectively: by what has been said this is
an acceptable abuse of notation. The same considerations
apply to the map $\xi:J \rightarrow A$, and the induced maps on
$p$-adic Tate modules (tensored with $W(k)$ and localised at $\sf m$) and the character
groups of tori (tensored with $W(k)$ and localised at $\sf m$) associated to the
special fiber at $q'$ of the N\'eron models of $J$ and $A$. 

We have the important Hecke equivariant exact sequence
of Ribet
\begin{eqnarray}\label{ribet}
0 \rightarrow {\cal X}(J,q')_{{\sf m},W(k)}
\rightarrow {\cal X}(J',q)_{{\sf m},W(k)} \rightarrow {\cal X}(J'',q)_{{\sf m},W(k)}^2
\rightarrow 0
\end{eqnarray} 
where ${\cal X}(J'',q)_{{\sf m},W(k)}$ is the character group of the torus
associated to the fibre at $q$ (tensored with
$W(k)$ and localised at $\sf m$) of the N\'eron model of the Jacobian of the Shimura
curve arising from the congruence subgroup
$\Gamma_0(N(\rhobar)p^{\delta}q_{2(r-s)} \cdots
q_{2r}/t)$ of $B_{q_1,\cdots,q_{2(r-s)-2}}$. We will denote the map ${\cal X}(J,q')_{{\sf m},W(k)}
\rightarrow {\cal X}(J',q)_{{\sf m},W(k)}$ of the exact sequence
(\ref{ribet}) by $\iota$.
The exact sequence (\ref{ribet}) is deduced easily from Proposition 1 of [RT]
(see also the original source, i.e., Theorem 4.1 of [Ri]). 
Let us denote by the same symol ${\cal L}$ the submodules 
${\cal X}(J)_{{\sf m},W(k)}[{\rm ker}(\pi)]$ and ${\cal
X}(J')_{{\sf m},W(k)}[{\rm ker}(\pi)]$ (free of rank one over $W(k)$) 
of ${\cal X}(J)_{{\sf m},W(k)}$ and ${\cal
X}(J')_{{\sf m},W(k)}$ respectively. 
This is an acceptable abuse of notation because the torsion-freeness of the
cokernel of $\iota$, that follows from (\ref{ribet}), implies that that these
different ${\cal L}$'s are identified by $\iota$. 

Recall the $W(k)$-bilinear monodromy pairings (see [RT]) 
$$(\ , \ )_J:{\cal X}(J,q')_{{\sf m},W(k)} \times {\cal X}(J,q')_{{\sf m},W(k)} \rightarrow
W(k)$$
$$(\ , \ )_{J'}:{\cal X}(J',q)_{{\sf m},W(k)} \times {\cal X}(J',q)_{{\sf m},W(k)} \rightarrow
W(k)$$ for which the action  of the operators $T_r$ is self-adjoint, and 
that are compatible with the injection
$$\iota: {\cal X}(J,q')_{{\sf m},W(k)} \rightarrow {\cal
X}(J',q)_{{\sf m},W(k)}, $$ 
of the exact sequence (\ref{ribet}) above. We also have the $W(k)$-bilinear monodromy pairings
$$(\ , \ )_{A'}:{\cal X}(A',q)_{{\sf m},W(k)}  \times {\cal X}(A',q)_{{\sf m},W(k)}
\rightarrow W(k)$$
$$(\ , \ )_A:{\cal X}(A,q')_{{\sf m},W(k)} \times{\cal X}(A,q')_{{\sf m},W(k)}
\rightarrow W(k)$$ where we are using the ismorphisms between ${\cal
X}(A',q)_{{\sf m},W(k)}$ and ${\cal X}((A')^d,q)_{{\sf m},W(k)}$ and ${\cal
X}(A,q')_{{\sf m},W(k)}$ and ${\cal X}((A)^d,q')_{{\sf m},W(k)}$
discussed earlier.  
We denote the torsion cokernels of the corresponding maps into
the duals, i.e., the component groups, 
by $\phi_{q'}(J)_{{\sf m},W(k)}$, $\phi_{q}(J')_{{\sf m},W(k)}$, 
$\phi_{q'}(A)_{{\sf m},W(k)}$ and $\phi_{q}(A')_{{\sf m},W(k)}$. 
(We can define the component groups at primes in $S$ analogously.)
The maps $\xi$ and $\xi'$ also induce maps 
$\xi_*:\phi_{q'}(J)_{{\sf m},W(k)} 
\rightarrow \phi_{q'}(A)_{{\sf m},W(k)}$
and
$\xi'_*:\phi_{q}(J')_{{\sf m},W(k)} \rightarrow \phi_{q}(A')_{{\sf m},W(k)}$
that are considered in the following proposition.

Before stating the proposition note that by 
invoking the result of [Ca] we see that the $GL_2(W(k))$-valued
representation of $G_{\bf Q}$ that arises from ${\rm Ta}_p(A)_{{\sf m},W(k)}
\simeq ({\rm Ta}_p(A) \otimes W(k))_{\sf m}$ depends only
on the (Hecke equivariant, ${\bf Q}$-)isogeny class of $A$ and is isomorphic to
$\rho_Q^{\Qnew}$ ($\simeq {\rm Ta}_p(A')_{{\sf m},W(k)}$).

\begin{prop}\label{component}
  We have the following information about maps between component
  groups:
\begin{enumerate}
\item We have the equalities of orders of component groups:
           $|\phi_{q}(A)_{{\sf m},W(k)}|=|W(k)/(x_{q})|=|\phi_{q}(A')_{{\sf
           m},W(k)}|$($x_{q}$ was defined in 3.1)
           for any prime $q \in Q$. Further, for any prime
           $\ell|N(\rhobar)p^{\delta}$, 
           the component groups $\phi_{\ell}(A')_{{\sf m},W(k)}$ and 
           $\phi_{\ell}(A)_{{\sf m},W(k)}$ are 
           trivial. The orders of these component groups depend
           only on the (Hecke equivariant, ${\bf Q}$-)isogeny class
           of $A$ (or $A'$).

\item  The group $\phi_q(J')_{{\sf m},W(k)}$ is trivial
       and hence the map $\xi'_*:\phi_{q}(J')_{{\sf m},W(k)} \rightarrow 
       \phi_q(A')_{{\sf m},W(k)}$ is trivial.


\item The map $\xi_*:\phi_{q'}(J)_{{\sf m},W(k)} \rightarrow
           \phi_{q'}(A)_{{\sf m},W(k)}$ is surjective.

\end{enumerate}
\end{prop}

\noindent{\bf Proof:} 
The first part is well-known and follows easily from using non-archimedean uniformisations
of the abelian varieties $A$ and $A'$ that have purely toric
reduction at all $\ell \in Q $ and $\ell|N(\rhobar)p^{\delta}$ (see Chapter III
of [Ri3] and Section 3 of [GS] for instance). 
The referee has remarked that this also follows directly from
the definition of the component group and the N\'eron mapping property.
The point is that 
$\phi_q(A')[\lambda^n] \cong W(k)/p^n$ if and only if $A'[\lambda^n]$
extends to a finite flat group scheme over ${\bf Z}_q$ (in
the notation of the paragraph at the end of Section 3.1): see the proof of Lemma~6.2 of [Ri].
The last sentence of part 1 follows from all these considerations and the sentence
before the statement of the proposition.

That $\phi_q(J')_{{\sf m},W(k)}$ is trivial follows 
as $\sf m$ is non-Eisenstein (see Theorem 3.12 of [Ri] for
an argument for this in a different situation but that works in ours too).

From part 1 it follows that the component group $\phi_t(A)_{{\sf m},W(k)}$ 
is trivial. 
Thus part 3 follows from part 1 and from the methods of 
Section 3 of [RT] (see Proposition 3 and its corollary of loc. cit.)
which show that the order of the cokernel
of $\xi_*:\phi_{q'}(J)_{{\sf m},W(k)} \rightarrow \phi_{q'}(A)_{{\sf m},W(k)}$
divides the order of  $\phi_{q_1}(A)_{{\sf
m},W(k)}=1$ (note
that our chosen prime $t=q_1$ divides the discriminant of the the quaternion algebra
from which $J$ arises for all values of $s$ ($0 \leq s \leq r-1$)
that we consider which makes the appeal to loc. cit. valid).
This finishes the proof of the proposition.

\vspace{3mm}

Now arguing exactly as in Section 2 of [RT] we use 
the exact sequence (\ref{ribet}) and the Proposition 
above to deduce Theorem \ref{RT}. (Most parantheses below will
indicate ideals of $W(k)$.) Consider ${\cal L}={\cal
X}(J',q)_{{\sf m},W(k)}[\pi]$: choose a generator $l$ of the $W(k)$-module $\cal L$
and set $(\tau):=((l,l)_{J'})$ where we are using the monodromy
pairing recalled above. We have the maps as above $\xi'^*:{\cal X}(A',q)_{{\sf m},W(k)}
\rightarrow {\cal X}(J',q)_{{\sf m},W(k)}$
and  $\xi'_*:{\cal X}(J',q)_{{\sf m},W(k)} \rightarrow {\cal X}(A',q)_{{\sf m},W(k)}$
(the second map results from the isomorphism chosen earlier between
${\cal X}((A')^d,q)_{{\sf m},W(k)}$ and ${\cal X}(A',q)_{{\sf
m},W(k)}$). Now 
by the adjointness property of $\xi'$ with respect to the monodromy
pairings above, i.e., $((\xi'^*x,y)_{J'})=((x,\xi'_*y)_{A'})$ for
any $x \in {\cal X}(A',q)_{{\sf m},W(k)}, y \in {\cal X}(J',q)_{{\sf
m},W(k)}$ (see discussion before Proposition 2 of [RT]) we have that 
$(\xi'_*\xi'^*)((x,x)_{A'})=((\xi'^*x,\xi'^*x)_{J'})$. For a generator $x$
of the free $W(k)$ rank-one module ${\cal X}(A',q)_{{\sf m},W(k)}$, we have that $((x,x)_{A'})=(x_q)$ from
Proposition \ref{component}, and thus we deduce that
$$(\xi'_*\xi'^*)(x_q)=([{\cal L}:{\cal X}(A',q)_{{\sf
m},W(k)}]^{1/d})^2.(\tau)$$
where $d$ is the rank of $W(k)$ as a ${\bf Z}_p$-module, and where we
view ${\cal X}(A',q)_{{\sf m},W(k)}$ as embedded in ${\cal L}$ by $\xi'^*$.
Arguing similarly with respect to the map $\xi:J \rightarrow A$ and
this time using the character groups of the tori of the reduction of
$J$ and $A$ mod
$q'$ and the equality $((\iota^{-1}(l),\iota^{-1}(l))_{J})=(\tau)$, we get the analogous equality 
$$(\xi_*\xi^*)(x_{q'})=([{\cal L}:{\cal X}(A,q')_{{\sf m},W(k)}]^{1/d})^2.(\tau).$$ Eliminating
$\tau$ from these equalities we get the equation of fractional ideals
of $W(k)$, $${ {(\xi'_*\xi'^*)(x_q)} \over  {([{\cal L}:{\cal
X}(A',q)_{{\sf m},W(k)}]^{1/d})^2}      }  =  { {(\xi_*\xi^*)(x_{q'}) } \over 
{ ([{\cal L}:{\cal X}(A,q')_{{\sf m},W(k)}]^{1/d})^2 } }.$$ Theorem \ref{RT} will follow from this
and Proposition \ref{component} once we
have proved the following lemma (see Proposition 2 of [RT]).

\begin{lemma}
 The index $[{\cal L}:{\cal X}(A',q)_{{\sf m},W(k)}]$ 
 is the cardinality of the cokernel of $\xi'_*:\phi_{q}(J')_{{\sf m},W(k)} \rightarrow 
 \phi_q(A')_{{\sf m},W(k)}$ and the index $[{\cal L}:{\cal
 X}(A,q')_{{\sf m},W(k)}]$ 
  is the cardinality of the cokernel of $\xi_*:\phi_{q'}(J)_{{\sf m},W(k)} \rightarrow 
 \phi_{q'}(A)_{{\sf m},W(k)}$.
\end{lemma}

\noindent{\bf Proof:} The optimality of $\xi':J' \rightarrow A'$ 
implies that the dual map $(A')^d \rightarrow J'$ is injective
and thus the map of character groups $\xi'_*:{\cal X}(J',q)_{{\sf m},W(k)}
\rightarrow {\cal X}((A')^d,q)_{{\sf m},W(k)}$ is surjective. We use the isomorphism
between ${\cal X}(A',q)_{{\sf m},W(k)}$ and ${\cal X}((A')^d,q)_{{\sf m},W(k)}$ to regard $\xi'_*$
as before as a map $\xi'_*:{\cal X}(J',q)_{{\sf m},W(k)}
\rightarrow {\cal X}(A',q)_{{\sf m},W(k)}$. Again using this isomorphism
and the prior description of the monodromy pairings we have
the commutative diagram:

\begin{equation}
\begin{CD}
0 @> >>{\cal X}(J',q)_{{\sf m},W(k)} @> >> {\rm Hom}_{W(k)}({\cal X}(J',q)_{{\sf m},W(k)},W(k))
 @> >> \phi_{q}(J')_{{\sf m},W(k)}  @>  >> 0
\\ 
&& @V VV @V VV @V VV 
\\
0 @> >>{\cal X}(A',q)_{{\sf m},W(k)} @> >> {\rm Hom}_{W(k)}({\cal X}(A',q)_{{\sf m},W(k)},W(k))
 @> >> \phi_{q}(A')_{{\sf m},W(k)}  @>  >> 0
\end{CD}
\end{equation} 
The first vertical map (which is the $\xi'_*$ above) being
surjective,
the cokernels of ${\rm Hom}_{W(k)}(\xi'^*,W(k))$ (which is the middle
vertical map) and the right hand vertical map (which is the map 
$\xi'_*:\phi_{q}(J')_{{\sf m},W(k)} \rightarrow 
 \phi_q(A')_{{\sf m},W(k)}$ we are interested in) may be identified.
But now the order of the cokernel of ${\rm Hom}_{W(k)}(\xi'^*,W(k))$
coincides with the order of the torsion subgroup of the cokernel of
the map $\xi'^*:{\cal X}(A',q)_{{\sf m},W(k)} \rightarrow {\cal X}(J',q)_{{\sf m},W(k)}$.
But since ${\cal X}(J',q)_{{\sf m},W(k)}/{\cal L}$ is torsion-free we derive the equality
of the index $[{\cal L}:{\cal X}(A',q)_{{\sf m},W(k)}]$ 
with the cardinality of the cokernel of $\xi'_*:\phi_{q}(J')_{{\sf m},W(k)} \rightarrow 
\phi_q(A')_{{\sf m},W(k)}$ as desired. The second statement of the lemma is proved
entirely analogously and this finishes the proof of the lemma.

\subsection{End of proof}

We start with the isomorphism 
$R_Q^{\Qnew} \simeq {\bf T}_Q^{\Qnew}$ of Theorem 1.
We apply Theorem 4 inductively to get the containment $(\xi_{r,*}\xi_r^*)
\subset (\prod_{q\in Q} x_q)$ of ideals of $W(k)$. 
Consider the monodoromy pairing on
$X = {\cal X}(J_r,t)_{{\sf m},W(k)}$ that was considered earlier where
$J_r=J_0(N(\rhobar)p^{\delta}q_2\cdots q_{2r})$: 
this is perfect as $\phi_{t}(J_r)_{{\sf m},W(k)}=0$ (see
part 2 of Proposition 3).
By Lemma 2 and part 1 of Prop.~3, $\xi_r^*$ defines a map from
${\cal X}(A_r,t)_{{\sf m},W(k)}$
to $X$ with torsion-free cokernel. 
Thus $\xi_r^*({\cal X}(A_r,t)_{{\sf m},W(k)})={\cal L}:=X[{\rm ker}(\pi)]$. 
Picking a generator $x$ of the free rank one $W(k)$-module
${\cal X}(A_r,t)_{{\sf m},W(k)}$, we see that $\xi_r^*(x)$
generates ${\cal L}$, and that
$((\xi_r^*(x),\xi_r^*(x)))=((\xi_{r,*}\xi_r^*)(x,x))=((\xi_{r,*}\xi_r^*))$
(equality of ideals of $W(k)$) 
where the pairings are the monodromy pairings on $X$ and ${\cal
  X}(A_r,t)_{{\sf m},W(k)}$ respectively.
The last equality follows from the fact that the component
group $\phi_t(A_r)_{{\sf m},W(k)}$ is trivial by part 1 of Proposition 3. 
Now just as in proof of Lemma 4.17 of [DDT]
we know that 
$X/(X[\ker(\pi)] + X[I])$ is a $W(k)/\pi(I)$-module
(where $I = {\rm Ann}_{{\bf T}_Q}({\rm ker}(\pi))$)  isomorphic
to the cokernel of ${\cal L} \rightarrow {\rm Hom}_{W(k)}({\cal
L},W(k))$ (the map being induced by the monodromy pairing on $X$)
since $X/X[I] \cong {\rm Hom}_{W(k)}(X[{\rm ker}(\pi)],W(k))$.
We conclude that $(\pi({\rm Ann}_{{\bf T}_{Q}}({\rm ker}(\pi)))) \subset (\prod_{q\in Q} x_q)$ 
and thus by Proposition 2 and  Theorem 3 we deduce the isomorphism 
$R_Q \simeq {\bf T}_Q$.

\vspace{3mm}

\noindent{\bf Remarks.}

\begin{itemize}

\item The presentation in 3.3 above using pairings on character groups
was suggested by the referee. We had earlier used pairings on Tate modules.
But note that to prove all the 
intermediate isomorphisms $R_{Q}^{Q_s-{\rm new}} \rightarrow {\bf
T}_{Q}^{Q_s-{\rm new}}$ (with $Q_s=\{q_1,\cdots,q_{2(r-s)}\}$) from the methods here one would in 3.3 
above have to use pairings on Tate module of $A_s$. The pairings
on character groups would not serve the purpose 
as the monodromy pairings would not be perfect (they would also not have 
faithful actions of the Hecke algebras that would 
be relevant unless $s=r$). 
On the other hand all these intermediate isomorphisms can be deduced
from $R_Q \simeq {\bf T}_Q$ using the methods of the appendix [B].

\item Note that in the proof of Theorem \ref{main} we did 
  not use any {\it a priori} multiplicity one
 results or Gorenstein properties of the Hecke algebras
 considered here. In [D] the construction in [TW] too was made independent of this
 input.  

\item In the proof the assumption $p>5$ is needed
only as far as this is needed in [R1]. It is quite likely that
the method of proof of Theorem
\ref{main} here  be extended
to $p>2$ (using the refinements of [R1] in [T]).

\item Richard Taylor had sketched to the author a monodromy argument for
computing change of $\eta$-invariants when one drops one prime at a
time that was very helpful in showing the way towards the proof given above.

\end{itemize}

\section{Level substitution mod $p^n$}

We briefly sketch in this section another approach
to lifting an isomorphism between new quotients of deformation and Hecke
rings that works under the assumption
that $N(\rhobar)p^{\delta}$ is divisible by at least 2 primes
(required for D. Helm's result in the appendix [H]).
The {\it raison d'etre} of this section is
to indicate a method of ``level substitution mod $p^n$'' 
(see Proposition \ref{sub}) and we will skip details that can be filled in by 
standard arguments.

For any finite set of primes $Q=\{q_1,\cdots, q_n\}$ 
that are special for $\rhobar$ (see Definition \ref{special})
and a subset $\alpha$ of $Q$
consider $H^1(X_0(N(\rhobar)p^{\delta}Q),W(k))^{\alpha-{\rm new}}$
which is defined as the maximal torsion-free quotient of the quotient of
$H^1(X_0(N(\rhobar)p^{\delta}Q),W(k))$
by the $W(k)$-submodule spanned by the images of 
$H^1(X_0(N(\rhobar)p^{\delta}{Q 
\over q}),W(k))^{2}$ in $H^1(X_0(N(\rhobar)p^{\delta}Q),W(k))$, 
as $q$ runs through the primes of $\alpha$, 
under the standard degeneracy maps (here and below for a finite set of primes $Q$
we abusively denote by $Q$ again the product of the primes in it). 
We consider the standard action of Hecke operators $T_r$ for 
all primes $T_r$ (note that we are using $T_r$ for operators 
that sometimes 
get called $U_r$ to be consistent with Helm's appendix). By [DT] there is a maximal ideal $\sf m$
of the $W(k)$-algebra generated by the action of these $T_r$'s such that
$T_r -\widetilde{{\rm tr}(\rhobar({\rm
Frob}_r))} \in \sf m$ for $\widetilde{{\rm tr}(\rhobar({\rm
Frob}_r))}$ an arbitrary lift to $W(k)$ of ${\rm tr}(\rhobar({\rm
Frob}_r))$ and $(r,N(\rhobar)pQ)=1$, and $T_r -\tilde{\alpha_r} \in \sf m$
for $r|N(\rhobar)pQ$, where $\alpha_r$ is the unique root
of the characteristic polynomial of $\rhobar({\rm Frob}_r)$ congruent to $\pm 1$ when 
$r|Q$, it is the scalar by which (the arithmetic Frobenius) 
${\rm Frob}_r$ acts on the unramified quotient
of $\rhobar|_{G_r}$ when $r|N(\rhobar)$ or $r=p$ is ordinary for $\rhobar$
(if $r=p$ and $\rhobar$ is not ordinary at $p$
we take $\alpha_r$ to be 0), and $\tilde{\alpha_r}$ is any lift
of $\alpha_r$ to $W(k)$. Then we define ${\bf T}_Q^{\alpha-{\rm new}}$ to be the localisation at $\sf m$
of the $W(k)$-algebra generated by the action of these Hecke operators on
the finite flat $W(k)$-module $H^1(X_0(N(\rhobar)p^{\delta}Q),W(k))^{\alpha-{\rm new}}$.
An analog of Lemma 1 gives that we have natural surjective maps
$R_Q^{\alpha-{\rm new}} \rightarrow {\bf T}_Q^{\alpha-{\rm new}}$
(where we take care of the fact that {\it all} 
$T_r$'s are in the image, including $r|N(\rhobar)p^{\delta}Q$, 
as in Section 2 of [W] relying on results of [Ca1]).

We need properties of ${\bf T}_Q^{\alpha-{\rm new}}$
that are proved ([H] and Corollary \ref{DT1}) 
by exploiting another description of these
algebras that we recall for orienting the reader (although we do not make explicit use 
of the alternative descriptions in this sketch). For this fix a subset $\beta$ of the prime divisors of
$N(\rhobar)p^{\delta}$. Denote by $B_{\alpha,\beta}$ the quaternion algebra over ${\bf Q}$ 
ramified at the primes in $\alpha \cup \beta$ and further at $\infty$
if the cardinality $n'=|\alpha \cup \beta|$ is odd. 
Denote by ${\bf A}$ the adeles over ${\bf Q}$.
For the standard open compact subgroup 
$U_{\alpha,\beta}:=U_0(N(\rhobar)p^{\delta}Q\alpha^{-1}\beta^{-1})$ 
of the ${\bf A}$-valued points of the algebraic group 
$G_{\alpha,\beta}$ (over ${\bf Q}$) corresponding to $B_{\alpha,\beta}$, $G_{\alpha,\beta}({\bf A})$, we consider the coset 
space ${\cal X}_{U_{\alpha,\beta}}=G_{\alpha,\beta}({\bf Q})\backslash G_{\alpha,\beta}({\bf A}) /U_{\alpha,\beta}.$ 
Depending on 
whether $n'$ is odd or even, this double coset space either
is merely a finite set of points, or can be given the structure of 
a Riemann surface (that is compact if $n' \neq 0$ and can be 
compactified by adding finitely many points if $n'$ is 0).
If $n'$ is odd we consider the space of functions  ${\cal 
S}_{U_{\alpha,\beta}}:=$
$\{f:{\cal X}_{U_{\alpha,\beta}} \rightarrow 
W(k)\}$  
modulo the functions which factorise through the norm map, 
and in the case of $n'$ even we consider 
the first cohomology of the corresponding Riemann surface ${\cal X}_{U_{\alpha,\beta}}$,
i.e., ${\cal S}_{U_{\alpha,\beta}}:=H^1({\cal X}_{U_{\alpha,\beta}},W(k))$. These  
$W(k)$-modules have the standard action of Hecke operators 
$T_r$ (see Sections 2 and 3 of [DT1]). By the results
of [DT] and the Jacquet-Langlands correspondence there is
a maximal ideal that we denote by $\sf m$ again in the support of the $W(k)$-algebra
generated by the action of the $T_r$'s on ${\cal S}_{U_{\alpha,\beta}}$ characterised as before. 
We denote the localisation at $\sf m$ of this Hecke
algebra by ${\bf T'}_Q^{\alpha \cup \beta-{\rm new}}$.
Then by the Jacquet-Langlands correspondence, which gives an
isomorphism ${\bf T'}_Q^{\alpha \cup \beta-{\rm new}} \otimes {\bf Q}_p \simeq {\bf
T}_Q^{\alpha-{\rm new}} \otimes {\bf Q}_p$ that takes $T_r$ to $T_r$, and the freeness of
${\bf T'}_Q^{\alpha \cup \beta-{\rm new}}$ and ${\bf T}_Q^{\alpha-{\rm
new}}$ as $W(k)$-modules, we have 
${\bf T'}_Q^{\alpha \cup \beta-{\rm new}} \simeq {\bf
T}_Q^{\alpha-{\rm new}}$, an isomorphism of local $W(k)$-algebras.

We consider an auxiliary set $Q=\{q_1,\cdots,q_n\}$ as in Theorem 1 such that
$R_Q^{\Qnew} \simeq {\bf T}_Q^{\Qnew} \simeq W(k)$, denote by $\pi'$ all the asserted isomorphisms, 
and define the ideal $(x_{q_i})$ of $W(k)$  as before in Section 3.1. We may
assume that  $Q \backslash \{q_i\}$
is not an auxiliary set (otherwise work with $Q \backslash
\{q_i\}$). Then by the
methods of [R], [R1], [T] and [KR] we have the following proposition.

\begin{prop}\label{sub}
  For each $q_i \in Q$ such that $x_{q_i} \in p^{m_i}W(k)
  -p^{m_i+1}W(k)$ there is a prime $q_i'$ not in $S \cup Q$ which is special
  for $\rho_Q^{\Qnew}$ mod $p^{m_i}$ and such that $Q_i=Q \backslash
  \{q_i\} \cup \{q_i'\}$ is an auxiliary set. Hence we have 
  isomorphisms 
  $R_{Q_i}^{Q_i-{\rm new}} \simeq {\bf T}_{Q_i}^{Q_i-{\rm new}} \simeq
  W(k)$  and further $\rho_Q^{\Qnew} \simeq \rho_{Q_i}^{Q_i-{\rm new}}$ mod
  $p^{m_i}$.
\end{prop}

\noindent{\bf Proof:} We use Lemma 8 of the companion paper [KR]. 
As $Q \backslash \{q_i\}$ is not an
auxiliary set, it follows using the notation of
loc. cit. that $H^1_{{\cal N}_v}(G_{S \cup Q_i},{\rm Ad}^0(\rhobar))$
and $H^1_{{\cal N}_v^{\perp}}(G_{S \cup Q_i},{\rm Ad}^0(\rhobar)(1))$
are 1-dimensional spanned by $\psi$ and $\phi$ respectively. 
Then using loc. cit. we 
choose $q_i'$ not in $S \cup Q$ 
such that $\psi|_{G_{q_i'}} \neq 0$, $\phi|_{G_{q_i'}} \neq 0$,
and $\rho_Q^{\Qnew}$ mod $p^{m_i}$ is special at $q_i'$.
This is the $q_i'$ that we want.

\vspace{3mm}

Set $Q'=\{q_1',\cdots,q_n'\}$. From the result of [H] and the isomorphisms
${\bf T'}_Q^{\alpha \cup \beta-{\rm new}} \simeq {\bf
T}_Q^{\alpha-{\rm new}}$ it follows that
the Hecke algebras ${\bf T}_{Q \cup Q'}^{\alpha-{\rm new}}$ are Gorenstein
for any subset $\alpha$ of $Q$. 
(In [H], Gorensteinness results are proved for the Hecke algebras that arise from action
on certain Tate modules: these imply (using [DT]) Gorensteinness of the Hecke
algebras considered here.
For the result of [H] the fact that $\rhobar({\rm Frob}_{q_i})$ is
not a scalar for $q_i \in Q$ is crucial, and the 
technical assumption that $\rhobar$ is not finite at
(at least) 2 primes is needed.)
Using \begin{itemize} \item level-raising methods of [DT] and [DT1] (that prove the analog of the Ihara-Ribet lemma
in the setting of Shimura curves) and again invoking the isomorphism ${\bf T'}_Q^{\alpha \cup \beta-{\rm new}} \simeq {\bf
T}_Q^{\alpha-{\rm new}}$, \item a simple Galois cohomology computation (see Section 4.2 of [DR]) 
\item and the numerical isomorphism criterion of
Wiles (Theorem \ref{numerical criterion}) \end{itemize} 
it is quite standard (after Section 2 of [W])
to deduce the following result
from the isomorphism $\pi':R_Q^{\Qnew} \simeq {\bf T}_Q^{\Qnew}$.
We skip its proof in this sketch.

\begin{cor}\label{DT1}
  We have an isomorphism $R_{Q \cup Q'}^{\Qnew} \simeq {\bf T}_{Q \cup
  Q'}^{\Qnew}$.
\end{cor}

For any subset $\alpha$ of
$Q$, we denote by $\pi$ the morphism ${\bf T}_{Q \cup
Q'}^{\alpha-{\rm new}} \rightarrow W(k)$ that sends $T_r$ to $\pi'(T_r)$
for $r \notin Q'$, but sends $T_r$ to the
unique eigenvalue of $\rho_Q^{\Qnew}({\rm Frob}_r)$ that is congruent to either plus or minus 1
for $r \in Q'$. For any subset $\alpha$ of $Q$ that does not contain $q_i$,
we denote by $\pi_i$ the morphism ${\bf T}_{Q \cup
Q'}^{\alpha-{\rm new}} \rightarrow W(k)$, that arises
from the isomorphism $\pi_i':{\bf T}_{Q_i}^{{Q_i}-{\rm new}} \simeq W(k)$ of Proposition \ref{sub}, that sends $T_r$ to $\pi_i'(T_r)$
for $r \notin Q' \cup \{q_i\} \backslash \{q_i'\}$ and sends $T_r$ to the
unique eigenvalue of $\rho_{Q_i}^{Q_i-{\rm new}}({\rm Frob}_r)$ that is congruent to either plus or minus 1
for $r \in Q' \cup \{q_i\} \backslash \{q_i'\}$.

Starting from Corollary \ref{DT1}, and using the morphisms $\pi_i$ for $i=1,\cdots,n$,
Theorem \ref{numerical criterion}, Proposition \ref{tangent space},
we briefly indicate how we deduce the isomorphism $R_{Q \cup Q'} \simeq {\bf T}_{Q \cup
Q'}$. To do this inductively it suffices to prove that
$R_{Q \cup Q'}^{Q_s-{\rm new}} \simeq {\bf T}_{Q \cup Q'}^{Q_s-{\rm new}}$
implies $R_{Q \cup Q'}^{Q_{s+1}-{\rm new}} 
\simeq {\bf T}_{Q \cup Q'}^{Q_{s+1}-{\rm new}}$ where 
$Q_s=\{q_1,\cdots,q_{n-s}\}$ and $s$ ranges over 0 to $n-1$.
As in the previous section to do this 
it will suffice to prove the following claim:

\vspace{3mm}

\noindent{\bf Claim:} The ideal 
$(\pi({\rm Ann}_{{\bf T}_{Q \cup Q'}^{Q_{s+1}-{\rm new}}}({\rm
  ker}(\pi))))$ of $W(k)$ is contained
in  $$(p^{m_{n-s}})(\pi({\rm Ann}_{{\bf T}_{Q \cup Q'}^{Q_{s}-{\rm
    new}}}({\rm ker}(\pi)))).$$

\vspace{3mm}

\noindent To prove the claim 
the two ingredients are the morphism $\pi_{n-s}:{\bf T}_{Q \cup Q'}^{Q_{s+1}-{\rm
    new}} \rightarrow W(k)$ which {\it does not} factor through ${\bf T}_{Q \cup Q'}^{Q_{s}-{\rm
    new}}$ by construction, and the result of [H] that
   ${\bf T}_{Q \cup Q'}^{Q_{s}-{\rm
    new}}$ and ${\bf T}_{Q \cup Q'}^{Q_{s+1}-{\rm
    new}}$ are Gorenstein.  
 For ease of notation set 
   ${\bf T}'={\bf T}_{Q \cup Q'}^{Q_{s+1}-{\rm new}}$
   and ${\bf T}={\bf T}_{Q \cup Q'}^{Q_{s}-{\rm new}}$ and denote the
   natural map ${\bf T}' \rightarrow {\bf T}$ by $\beta$. 

We justify the claim using these 2 ingredients: we
owe the following argument to the referee.
Because ${\bf T}$ and ${\bf T}'$ Gorenstein, we have that
$$\pi\beta({\rm Ann}_{{\bf T}'}({\rm ker}(\pi\beta)))=\pi\beta({\rm
Ann}_{{\bf T}'}({\rm ker}(\beta)))\pi({\rm Ann}_{{\bf T}}({\rm ker}(\pi)))$$ as
ideals of $W(k)$. Choose a $x \in {\rm ker}(\beta)$ such
that $\pi_{n-s}(x) \neq 0$ (this exists as $\pi_{n-s}$ 
does not factor through ${\bf T}$).
Consider $y \in {\rm Ann}_{{\bf T}'}({\rm ker}(\beta))$. Then as $xy=0$ this
implies $\pi_{n-s}(xy)=0$, which implies $\pi_{n-s}(y)=0$, 
which implies $\pi\beta(y) \in p^{m_{n-s}}W(k)$ as $\pi_{n-s}$ and $\pi$
are congruent mod $p^{m_{n-s}}$ which proves the claim.

\vspace{3mm}





\noindent{\bf Remarks:} \begin{itemize} \item It follows from [B] that
the isomorphism of complete intersections 
$R_Q \simeq {\bf T}_Q$ implies the isomorphism of complete intersections
$R_Q^{\alpha-{\rm new}} \simeq {\bf T}_Q^{\alpha-{\rm new}}$
for each subset $\alpha$ of $Q$ and hence in particular that
${\bf T}_Q^{\alpha-{\rm new}}$ is Gorenstein. Thus the Gorenstein
property of new quotients is a consequences of $R_Q \simeq {\bf T}_Q$-type theorems
which in turn can be proved (see the previous section) without 
knowing {\it a priori} such a property. In the method of
this section on the other hand we needed to use the Gorenstein property
of new quotients to prove $R_Q \simeq {\bf T}_Q$-type theorems.
\item The referee had pointed out that the arguments of this section in an
earlier version had a
 gap as they invoked  {\it a priori} Gorenstein properties of
 new-quotients of Hecke algebras that had not been justified. 
The result of [H] fills this gap.
\end{itemize}

\section{Acknowledgements} The author would like to express his
gratitude to Ravi Ramakrishna
for making available [R] and [R1] while still in preprint form and
for pointing out the striking consequence (Theorem \ref{DT} of this paper) 
of the method of [R1] that was the main inspiration for this work. 
He would also like to thank Fred Diamond and Richard Taylor for
helpful conversations, and Ken Ribet
for helpful correspondence. He would like to thank Gebhard B\"ockle and  David Helm
for writing appendices to the paper. This work germinated during a visit
to Brandeis University in the Fall of 1999 and he would like to thank Fred Diamond for
the invitation. He would like to thank the referee for a careful
reading and helpful critique of the article. 
Finally he would like to acknowledge the debt this article
owes to the beautiful ideas of Ken Ribet and Andrew Wiles.

\section{References}

\noindent [B] B\"ockle, G., {\it On the isomorphism $R_\emptyset \to T_\emptyset$}, 
appendix to this paper.

\vspace{3mm}

\noindent [Ca] Carayol, H., {\it Formes modulaires et repr\'esentations 
galoisiennes avec valeurs
dans un anneau local complet}, in {\it $p$-adic monodromy 
and the Birch and Swinnerton-Dyer conjecture} (Boston,
MA, 1991), 213--237, Contemp. Math., 165, AMS, 1994. 

\vspace{3mm}

\noindent [Ca1] Carayol, H., {\it Sur les repr\'esentations 
$\ell$-adiques associ\'ees aux
formes modulaires de Hilbert}, Annales de l'Ecole Norm. Sup.
19 (1986), 409--468.

\vspace{3mm}

\noindent [deS] de Shalit, E., {\it Deformation rings and Hecke rings},
in [FLT].

\vspace{3mm}

\noindent [D] Diamond, F., {\it The Taylor-Wiles 
construction and multiplicity
one}, Invent. Math. 128 (1997), no. 2, 379--391.

\vspace{3mm}

\noindent [DDT] Darmon, H., Diamond, F., Taylor, R., {\it 
Fermat's last theorem}, Current developments in mathematics, 1995 
(Cambridge, MA), 1--154, Internat. Press, Cambridge,
MA, 1994. 

\vspace{3mm}

\noindent [DR] Diamond, F., Ribet, K., {\it $\ell$-adic modular deformations
and Wiles's ``Main Conjecture''}, in [FLT].

\vspace{3mm}

\noindent [DT] Diamond, F., Taylor, R., 
{\it Lifting modular mod $l$ representations},
Duke Math. J. 74 (1994), no. 2, 253--269.

\vspace{3mm}

\noindent [DT1] Diamond, F., Taylor, R., {\it Nonoptimal 
levels for mod $\ell$ modular representations of 
$Gal(\overline{{\bf Q}}/{\bf Q})$}, Invent. Math. 115 (1994), 435--462.

\vspace{3mm}

\noindent [FLT] {\it Modular forms and Fermat's last theorem}, 
edited by Gary Cornell, Joseph H. Silverman and Glenn Stevens.
Springer-Verlag, New York, 1997.

\vspace{3mm}

\noindent [G] Grothendieck, A., Expos\'e IX, SGA7, SLNM 288, 1972. 

\vspace{3mm}

\noindent [GS] Greenberg, R., Stevens, G., {\it $p$-adic
$L$-functions 
and $p$-adic periods of modular forms},
Invent. Math. 111 (1993), no. 2, 407--447. 

\vspace{3mm}

\noindent [H] Helm, D., {\it The Gorenstein property for new quotients}, 
appendix to this paper.

\vspace{3mm}

\noindent [KR] Khare, C., Ramakrishna, R., {\it
Finiteness of Selmer groups and deformation rings}, preprint.

\vspace{3mm}

\noindent [Ma] Mazur, B., {\it An introduction to the deformation
theory of Galois representations}, in [FLT].

\vspace{3mm}

\noindent [R] Ramakrishna, R., {\it Lifting Galois
representations}, Invent. Math. 138 (1999), 537--562.

\vspace{3mm}

\noindent [R1] Ramakrishna, R., {\it Deforming Galois representations
and the conjectures of Serre and Fontaine-Mazur}, to appear in Annals
of Math.

\vspace{3mm}

\noindent [RT] Ribet, K., Takahashi, S., {\it 
Parametrizations of elliptic curves by
Shimura curves and by classical modular curves}, 
in {\it Elliptic curves and modular forms} (Washington, DC,
1996). Proc. Nat. Acad. Sci. U.S.A. 94 (1997), no. 21, 11110--11114.

\vspace{3mm}

\noindent [Ri] Ribet, K., {\it On modular representations of ${\rm
Gal}(\overline{\bf Q}/{\bf Q})$
arising from modular forms}, Invent. Math. 100 (1990), no. 2, 431--476.

\vspace{3mm}

\noindent [Ri1] Ribet, K., {\it Report on mod $\ell$
representations of ${\rm Gal}(\overline{\bf Q}/{\bf Q})$}, in Motives, Proc.
Sympos. Pure Math. 55, Part 2 (1994), 639--676.

\vspace{3mm}

\noindent [Ri2] Ribet, K., {\it Congruence relations between modular
forms}, Proc. International Cong. of Math. (1983), pp 503--514.

\vspace{3mm} 

\noindent [Ri3] Ribet, K., {\it Galois action on division
points of Abelian varieties with real multiplications},
Amer. J. Math. 98 (1976), no. 3, 751--804.

\vspace{3mm}

\noindent [Ri4] Ribet, K., {\it 
Multiplicities of Galois representations in Jacobians of Shimura
curves}, Israel Math. Conf. Proc., 3, Weizmann, Jerusalem, 1990.

\vspace{3mm}

\noindent [S] Serre, J-P., {\it  Sur les
repr\'esentations modulaires de degr\'e 2 de ${\rm
Gal}(\overline{\bf Q}/{\bf Q})$}, Duke Math. J. \textbf{54} (1987),
179--230.



\vspace{3mm}

\noindent[TW] Taylor, R., Wiles, A., 
{\it Ring-theoretic properties of certain Hecke
algebras}, Ann. of Math. (2) 141 (1995), no. 3, 553--572. 

\vspace{3mm}

\noindent [W] Wiles, A., {\it Modular elliptic curves and 
Fermat's last theorem}, Ann. of
Math. (2) 141 (1995), no. 3, 443--551. 

\vspace{5mm}        

\noindent {\bf Address of author:} School of Mathematics, 
TIFR, Homi Bhabha Road, Mumbai 400 005, INDIA: shekhar@math.tifr.res.in

\noindent Dept. of Math, University of Utah, 
155 S 1400 E, Salt Lake City, UT 84112, USA: shekhar@math.utah.edu

\end{document}